\documentclass{czjphys}
\usepackage{amsmath}
\usepackage{amsfonts}

\def\POZ{\footnote{Presented at the $13^{\rm th}$ International
Colloquium on ``Integrable Systems and Quantum Groups", Prague,
17--19 June 2004.}}

\newcommand {\tto} {\longrightarrow}
\usepackage{bbm}
\newcommand{\One}{\mathbbmss{1}}

\newcommand {\Cee}    {{\mathbb  C}}
\newcommand {\Zee}    {{\mathbb  Z}}
\newcommand {\Fee}    {{\mathbb  F}}

\newcommand {\fab}    {{\mathfrak{ab}}}
\newcommand {\fag}    {{\mathfrak{ag}}}

\newcommand {\fg}     {{\mathfrak{g}}}    %
\newcommand {\fh}     {{\mathfrak{h}}}

\newcommand {\fl}     {{\mathfrak{l}}}
\newcommand {\fL}     {{\mathfrak{L}}}

\newcommand {\fosp}   {{\mathfrak{osp}}}

\newcommand {\fpo}    {{\mathfrak{po}}}

\newcommand {\fpsl}   {{\mathfrak{psl}}}

\newcommand {\fsl}    {{\mathfrak{sl}}}

\newcommand {\fsvect} {{\mathfrak{svect}}}

\newcommand {\fvect}  {{\mathfrak{vect}}}   %

\newcommand {\bcdot}   {\mathbin{\hbox{\raise.4ex\hbox{\bf.}}}} 

\begin{document}

\title{Lie superalgebra structures in
\protect \bmth{H^{\bcdot}(\fg; \fg)}\,\POZ}
\authori{Pavel Grozman}
\addressi{Equa Simulation AB, Stockholm, Sweden; pavel@rixetele.com}
\authorii{Dimitry Leites}
\authoriii{}
\addressii{MPIMiS, Inselstr. 22, DE-04103 Leipzig, Germany\\
on leave from Department of Mathematics, University of Stockholm\\
Roslagsv.  101, Kr\"aftriket hus 6, SE-104 05, Stockholm, Sweden\\
mleites@math.su.se}
\headauthor{P. Grozman and D. Leites}
\headtitle{Lie superalgebra structures in $H^{\bcdot}(\fg; \fg)$}
\lastevenhead{P. Grozman and D. Leites: Lie superalgebra
structures in $H^{\bcdot}(\fg; \fg)$}

\refnum{A}
\daterec{23 July 2004}    
\issuenumber{11}  \year{2004} \setcounter{page}{1313}
\maketitle

\begin{abstract}
Let $\fg=\fvect(M)$ be the Lie (super)algebra of vector fields on
any connected (super)manifold $M$; let $\Pi$ be the change of
parity functor, $C^{i}$ and $H^{i}$ the space of $i$--chains and
$i$--cohomology.  The Nijenhuis bracket makes
$\fL_{\bcdot}=\Pi(C^{\bcdot-1}(\fg;\fg))=C^{\bcdot-1}(\fg;\Pi(\fg))$
into a Lie superalgebra that can be interpreted as the centralizer
of the exterior differential considered as a vector field on the
supermanifold $\hat M=(M, \Omega^{\bcdot}(M))$ associated with the
de Rham bundle on $M$.  A similar bracket introduces structures of
DG Lie superalgebra in $\fL_{\bcdot}$ and
$\fl_{\bcdot}=H^{\bcdot-1}(\fg; \Pi(\fg))$ for any Lie
superalgebra $\fg$.  We use a \textit{Mathematica}--based package
\textbf{SuperLie} (already proven useful in various problems) to
explicitly describe the algebras $\fl_{\bcdot}$ for some simple
finite dimensional Lie superalgebras $\fg$ and their ``relatives''
--- the nontrivial central extensions or derivation algebras of
the considered simple ones.
\end{abstract}

\section{Introduction}

This paper is a sequel to \cite{LLS} and \cite{GL3}. Its aim is to
demonstrate usefulness of the \textit{Mathematica}--based package
{\bf SuperLie} (\cite{SLie}) already tested on various problems
(\cite{GL2}). Our results --- computation of $H^*(\fg; \fg)$
(complete, or partial) --- are \lq\lq orthogonal" to those of
\cite{SZ} who computed only $H^1(\fg; M)$ and $H^2(\fg; M)$ but
for all (finite dimensional) irreducible modules $M$. We also
complete Tyutin's description \cite{Ty} of deformations of
$\fpo(0|n)$ who ignored odd parameters; cf. \cite{KST}.

During the talk at the conference we indicated how Lie algebra
cohomology is related to a \lq\lq new" (designed in late 1960's)
and seemingly never explored method for solving differential
equations and for the study of stability of various dynamical
systems. For an exposition of this and many other applications of
cohomology, see \cite{GL2}. For basics on Lie superalgebras we
consider, see \cite{LSc}.  For main data on Lie algebra cohomology
and rudiments of Lie superalgebra cohomology, see Fuchs's book
\cite{Fu}. It is well-known that if the Lie algebra $\fg$ is
finite dimensional and simple and $M$ a finite dimensional
$\fg$-module (both over $\Cee$), then
\begin{equation}
H^*(\fg; M)=0 \;\text{ for any irreducible $M$}\,. 
\end{equation}
However, there are very few general theorems helping to compute
$H^*(\fg; M)$ when both $\fg$ and $M$ are considered over $\Zee$
or $\Fee_q$, even less for Lie superalgebras.

On the other hand, C. Gruson suggested a totally new method
(applicable to Lie superalgebras only, not to Lie algebras)
\cite{Gr}.  At the moment, Gruson's method  is only applied to the
trivial coefficients. Its applicability in other cases is not
studied yet. So, for nontrivial modules over Lie superalgebras, we
have only the same tool researchers had at their disposal at the
birth of the cohomology theory: the definition. To see what
phenomena and patterns we might encounter, we use
\textbf{SuperLie} (\cite{SLie}) to get a supply of reliable
results to be used in more general analytic study.

The standard proof of $(1)$ uses the (even, quadratic) Casimir
operator. Passing to Lie superalgebras $\fg$  we observe that:
(A) some Lie superalgebras have no (even) quadratic Casimir operator,
(B) for some Lie superalgebras, such an operator exists but
vanishes at various $M$'s. The algebras we
consider here have these properties (A) and (B), and hence (1) fails sometimes.

In what follows, $C^i:=C^i(\fg; \fg)$ and $H^i:=H^i(\fg; \fg)$ for
a Lie superalgebra $\fg$; set $C^{\bcdot}=\oplus C^i$ and
$H^{\bcdot}=\oplus H^i$.

\section{The Lie superalgebra structure on \bmth{\Pi(C^{\bcdot-1})} and
\bmth{\Pi(H^{\bcdot-1})}}

Let the $e_\alpha$ be elements of the (weight) basis of $\fg$ (for
example, the Chevalley basis if $\fg$ is simple), and let the
$f^\alpha$ be the elements of the dual basis. The basis of
cochains is given by monomials of the form
$$
e_\alpha\otimes f^{\alpha_{1}}\wedge \ldots\wedge
f^{\alpha_{n}}\,.
$$
For $A=e^A_\alpha\otimes f_A^{\alpha_{1}}\wedge \ldots\wedge
f_A^{\alpha_{n_A}}\in C^{n_A}$ and $B=e^B_\beta\otimes f_B^{\beta_{1}}\wedge \ldots\wedge
f_B^{\beta_{n_B}}\in C^{n_B}$, we set
$$
A\cdot B=\sum_{k=1}^{n_{A}}
(-1)^{n_{A}-k}f_A^{\alpha_{k}}(e_{\beta}^B)e_{\alpha}^A\otimes
f_A^{\alpha_{1}}\wedge\dots\wedge \widehat{f_A^{\alpha_{k}}}\wedge
\dots\wedge f_A^{\alpha_{n_{A}}}\wedge f_B^{\beta_{n_{1}}}\wedge
\dots\wedge f_B^{\beta_{n_{B}}}\,.
$$
The {\it Nijenhuis bracket} on $\Pi(C^{\bcdot})$, where $\Pi$ is the shift
of parity functor, is given by the formula (hereafter $p(a)$ is the parity of
$a$ in $\Pi(C^{\bcdot})$, not in $C^{\bcdot}$)
\begin{equation}
[A, B]=A\cdot B-(-1)^{p(A)p(B)}B\cdot A\in C^{n_{A}+n_B-2}\,.
\end{equation}
Let us not only change parity, but also shift the degree by
setting $\deg A= n_{A}-1$; we denote this by writing
$\fL_{\bcdot}=\Pi(C^{\bcdot-1}):=\oplus_iC^{i-1}$. It is subject
to a direct verification that (2) defines the Lie superalgebra
structure on $\fL_{\bcdot}$, and
\begin{equation}
\D[A, B]=[\D A,B]+(-1)^{p(A)}[A,\D B]\,.
\end{equation}
Therefore $\Pi(Z^{\bcdot-1})$ is a subalgebra of
$\Pi(C^{\bcdot-1})$, and $\Pi(B^{\bcdot-1})$ is an ideal in
$\Pi(Z^{\bcdot-1})$.  Hence, we have a DG Lie superalgebras
structure on $\fL_{\bcdot}:=\Pi(C^{\bcdot-1})$ and
$\fl_{\bcdot}:=\Pi(H^{\bcdot-1})$. In \cite{LLS}, we observed that
the differential vanishes on $\fl_{\bcdot}$, and $\fl_{\bcdot}$ is
often very small, so although researchers mainly study
$\fl_{\bcdot}$, the algebra $\Pi(Z^{\bcdot-1})$ might be more
interesting than $\fl_{\bcdot}$ in some questions. Observe that
representing $\Pi(C^{\bcdot}(\fg;\fg))$ as $C^{\bcdot}(\fg;
\Pi(\fg))$ considerably simplifies computations.

\section{Examples}

In what follows, for small $i$, we listed superdimensions of the
$H^i$ expressed as $a|\bar b$; we write $a$ instead of $a|0$ and
$\bar b$ instead of $0|\bar b$. For $i$ greater than indicated
below, we did not calculate (partly due to
\textit{Mathematica}--imposed limitations).

For lists of simple Lie superalgebras, their $\Zee$--gradings,
known central extensions, outer derivations and deformations, see
\cite{K}, \cite{CK} and \cite{LSc}. On deformations of Poisson
superalgebras for various types of functions, see \cite{LSc},
\cite{Ty}, \cite{KST} and refs. therein. Here we only consider the
$0|n$--dimensional case. Recall that $\fvect(m|n)$ is the Lie
superalgebra of polynomial vector fields on the $m|n$--dimensional
superspace, $\fsvect(m|n)$ is its divergence--free subalgebra,
whereas $\fg=\fag_2$ and $\fab_3$ are exceptional finite
dimensional Lie superalgebras, most clearly determined by their
defining relations \cite{GL1}.

Every $\Zee$--grading of $\fg$ induces a natural $\Zee$--grading
on the space of all polynomial functions on $\fg$, in particular,
on cohomology. Such $\Zee$--grading will be called the
\textit{degree}. Observe that (1) the parity of a derivation
$c:\fg\tto\fg$ is opposite to the parity of $c$ considered as a
1--cocycle $c: \Pi(\fg)\tto\fg$; (2) to consider $\fl_{\bcdot}$,
one should shift \textit{all} parities in the tables below.
Observe that \textbf{the deformations with the odd parameter are
automatically global}.

\subsection{Simple Lie superalgebras without deformations}

\noindent
\textbf{Conjecture}. For $\fg=\fag_2$ and $\fab_3$, and
also $\fg=\fvect(0|n)$, where $n>1$, we have $\fl_{\bcdot}=0$.
(Verified to degree 5 and 4, and, for $n=2, 3$, to degree  7 and
5, respectively.)

\underline{$\fg=\fpsl(2|2):=\fsl(2|2)/\Cee \One_4$}.  This algebra
is highly symmetric: the group of its outer automorphisms is
Out~$\fg\simeq \fsl(2)$, so $\dim H^1= 3$.  We also knew that
$H^2=0$. Let $\deg$ be induced by the degree of $\fsl(2|2)$ in the
standard format, the even (diagonal) block matrices being of
degree 0, the upper/lower ones of degree $\pm 1$.  Here are known
and new results: listed are the superdimensions of the
$H^i=\fl_{i-1}$ for $i\leq 5$ and the degrees of their basis
elements: {\tiny
$$
\renewcommand{\arraystretch}{1.4}
\begin{tabular}{|c|c|c|c|c|c|c|c|}
\hline
$i-1\backslash \text{deg}$&$-6$&$-4$&$-2$&$0$&$2$&$4$&$6$\cr
\hline
$-1$& & & & & & &\cr
\hline
$0$& & &1&1&1&&\cr
\hline
$1$&&&&&&&\cr
\hline
$2$&&1&2&2&2&1&\cr
\hline
$3$&&&$\bar 1$&$\bar 1$&$\bar 1$&&\cr
\hline
$4$&1&2&2&2&2&2&1\cr
\hline
\end{tabular}
$$}

\noindent\textbf{Conjecture}. $\fl_{0}$, $\fl_{2}$, and $\fl_{3}$
generate $\fl_{\bcdot}$.

\subsection{Lie superalgebras with deformations}

\hskip1pc
\underline{$\fg=\fosp(4|2;\alpha)$}.  (Observe that, for
$\alpha=-1$ and 0, the algebra is not simple.)  The very
definition of the algebra as a 1--parameter family indicates that
$\dim H^2\geq 1$. Nothing was known about higher cohomology. Here
are new results: for $\alpha =1$ as well as for $\alpha =\frac37$
which served as generic $\alpha$ (the answer is the same, but for
the computer the task is much easier), we have
$$
\renewcommand{\arraystretch}{1.4}
\begin{tabular}{|c|c|c|c|c|c|c|c|}
\hline
$i-1$&$-1$&0&1&2&3&4&5\cr
\hline
$\dim H^{i-1}$&0&0&1&0&0&$\bar 2$&0\cr
\hline
\end{tabular}
$$

The cocycles depend on $\alpha$ but the product in $\fl_{\bcdot}$
does not (for  $\alpha\neq -1$, 0): the product of any two
cocycles from the above table vanishes.

\noindent \textbf{Conjecture}. There are no more cocycles but the
above.

\smallskip

The whole pathos of \cite{LLS} was to express $\fl_{\bcdot}$ as an
\textit{algebra}, not as a list of dimensions of the homogeneous
components, as classics did (Borel--Weil--Bott--Kostant). Here the
product of any two cocycles from the above table vanishes, so the
list is as fine as algebra as long as the product is trivial (even
if Conjecture fails). In examples that follow, we return, to our
regret, to the \lq\lq  list of dimensions" level.

\subsubsection{The divergence--free series}

We knew (\cite{LSc}) that $\fg=\fsvect(0|n)$ for $n>2$ has only
one outer derivation, so $\dim H^1= 1$, and $\fg$ has only one
global deformation, so $\dim H^2\geq 1$, but we knew nothing about
other cohomology. New results: indicated are nonzero $H^i$, in all
cases $H^i=\text{Span}(c_i)$ for some cocycles $c_i$, i.e., $\dim
H^i=1$ in all cases below. The parity of $c_i$, as an element of
$\fl_{\bcdot}$, is equal to $i-1+\text{deg}\pmod 2$.

\underline{$\fsvect(0|3)$}:
$$
\renewcommand{\arraystretch}{1.4}
\begin{tabular}{|c||c|c|c|c|c|}
\hline
$i-1$&$0$&$1$&$3$&$5$&$6$\cr
\hline
deg&$0$&$3$&$-3$&$0$&$-6$\cr
\hline
\end{tabular}
$$
Obviously, $[c_2, c_4]=0$, conjecturally $[c_4, c_4]=c_7$. With
\textit{Mathematica}'s inbuilt limitations, we were unable to
compute other products nor advance further than indicated (below
as well).

\underline{$\fsvect(0|4)$}:
$$
\renewcommand{\arraystretch}{1.4}
\begin{tabular}{|c||c|c|}
\hline
$i-1$&$0$&$1$\cr
\hline
deg&$0$&$4$\cr
\hline
\end{tabular}
$$

\subsubsection{The Poisson series $\fpo(0|m)$}

We knew (\cite{LSc}) that the infinite dimensional Lie
(super)algebra $\fg=\fh(2n|m)=\fpo(2n|m)/\text{center}$ of
Hamiltonian vector fields has only one outer derivation, and were
sure that the same holds for $\fg=\fpo(0|m)$; we also knew that
$\fg$ has only one global deformation, so $\dim H^2\geq 1$. We
knew nothing about other cohomology nor about the finite
dimensional case.

In order not to confuse the elements of $\fpo(0|m)$ with
functions, we realize $\fpo(0|m)$ as a subalgebra of contact
vector fields $K_f$ generated by functions $f$ (for formulas, see
\cite{LSc}). Let $\deg$ be the degree of an element of $\fg$ (and
the induced degree of cochains) given by the formula
$\deg(K_f)=\deg_{\mathrm{pol}}(f)-2$, where $\deg_{\mathrm{pol}}$
is the standard grading in the polynomial algebra (the degree of
each indeterminate is equal to 1). So, for $\fpo(0|n)$, there are
the following obvious cohomology: 0th, given by $K_1$; 1st, given
by the grading by the degree; 2nd, quantization, of $\deg=-4$;
and, additionally, the above cocycles wedged (which is possible
since $H^*(\fg; \fg)$ is an $H^*(\fg)$--module) by
$\left(d(K_{\theta_1 \cdot\ldots\cdot \theta_n})\right )^m$, where
$m$ is arbitrary for $n$ odd and is either 0 or 1 for $n$ even.
(Here $d(K_{\theta_1 \cdot\ldots\cdot \theta_n}): K_{\theta_1
\cdot\ldots\cdot \theta_n}\mapsto 1$.)

\underline{$\fg=\fpo(0|4)$}. {\tiny
$$
\renewcommand{\arraystretch}{1.4}
\begin{tabular}{|c|c|c|c|c|c|c|c|c|c|}
\hline
$i-1\backslash \deg$&$-8$&$-6$&$-4$&$-2$&$0$\cr
\hline
$-1$&&&&1&\cr
\hline
$0$&&&$\bar 1$&&$\bar 1$\cr
\hline
$1$&&&1&1&\cr
\hline
$2$&&$\bar 1$&&$\bar 3$&\cr
\hline
$3$&&1&3&&2\cr
\hline
$4$&$\bar 1$&&$\bar 3$&$\bar 2$&\cr
\hline
\end{tabular}
$$}

\underline{$\fg=\fpo(0|5)$}.  {\tiny
$$
\renewcommand{\arraystretch}{1.4}
\begin{tabular}{|c|c|c|c|c|c|c|c|c|c|c|c|c|c|c|c|c|}
\hline $i-1\backslash
\deg$&$-14$&$-13$&$-12$&$-11$&$-10$&$-9$&$-8$&$-7$&$-6$&$-5$&$-4$&$-3$&$-2$&$-1$&$0
$\cr \hline $-1$& & & & & & &&&&&&&1&&\cr \hline $0$& & & & &&&&
&&$1$&&&&&$\bar 1$\cr \hline $1$&&&&&&&1&&& &$1$&$\bar 1$&&&\cr
\hline $2$&&&&$1$&&&&$2$&$\bar 1$&&&&$\bar 2$&&\cr \hline
$3$&1&&&&2&$\bar 1$&&&1&$\bar 3$&&&&&1\cr \hline
\end{tabular}
$$}

\underline{$\fg=\fpo(0|6)$}.  {\tiny
$$
\renewcommand{\arraystretch}{1.4}
\begin{tabular}{|c|c|c|c|c|c|}
\hline
$i-1\backslash \deg$&$-8$&$-6$&$-4$&$-2$&$0$\cr
\hline
$-1$& & & &1&\cr
\hline
$0$& &$\bar 1$& & &$\bar 1$\cr
\hline
$1$&&&2&&\cr
\hline
$2$&$\bar 2$&&&$\bar 2$&\cr
\hline
\end{tabular}
$$}

\subsubsection{The Hamiltonian series $\fh(0|n)$ and
$\fh'(0|n)=[\fh(0|n), \fh(0|n)]$}

We knew (\cite{CK}) that the {\it infinite} dimensional
$\fg=\fh(2n|m)$  has only one outer derivation, so $\dim
H^1\geq1$, and (\cite{LSc}) that $\fg$ has only one global
deformation, except for $(2n, m)=(2, 2)$, so $\dim H^2\geq 1$, but
we knew nothing about other cohomology nor about the finite
dimensional case.  Let $\deg$ be the degree of an element of $\fg$
induced by the grading of $\fpo$. Here are new results:

\underline{$\fg=\fh(0|4)$}. {\tiny
$$
\renewcommand{\arraystretch}{1.4}
\begin{tabular}{|c|c|c|c|c|c|c|c|}
\hline
$i-1\backslash \deg$&$-6$&$-4$&$-2$&$0$&$2$&$4$\cr
\hline
$-1$& & & & & &  \cr
\hline
$0$& & & &$\bar 1$&&\cr
\hline
$1$&&1&&&&\cr
\hline
$2$&&&$\bar 1$&&$\bar 1$&\cr
\hline
$3$&1&1&&1&&\cr
\hline
$4$&&$\bar 2$&&&&$\bar 1$\cr
\hline
\end{tabular}
$$}

\underline{$\fg=\fh'(0|4)$}.  This is a special case already considered since
$\fh'(0|4)\simeq\fpsl(2|2)$.

\underline{$\fg=\fh(0|5)$}.  {\tiny
$$
\renewcommand{\arraystretch}{1.4}
\begin{tabular}{|c|c|c|c|c|c|c|c|c|c|c|c|c|c|}
\hline
$i-1\backslash
\deg$&$-10$&$-9$&$-8$&$-7$&$-6$&$-5$&$-4$&$-3$&$-2$&$-1$&$0$&1&$2$\cr
\hline
$-1$&&&&& & &&&&&&&\cr
\hline
$0$&&&&&&&& &&&$\bar 1$&& \cr
\hline
$1$&&&&& &&1&$\bar 1$&&&&&\cr
\hline
$2$&&&&1&$\bar 1$&&&&$\bar 1$&&& &$\bar 1$\cr
\hline
$3$&1&$\bar 1$&&&1&$\bar 1$&&& &$\bar 1$&&&\cr
\hline
\end{tabular}
$$}

\underline{$\fg=\fh'(0|5)$}.  {\tiny
$$
\renewcommand{\arraystretch}{1.4}
\begin{tabular}{|c|c|c|c|c|c|c|c|c|c|c|c|c|}
\hline
$i-1\backslash
\deg$&$-6$&$-5$&$-4$&$-3$&$-2$&$-1$&$0$&1&$2$&3&$4$&5\cr
\hline
$-1$& & &&&&&&&&&&\cr
\hline
$0$&&&& &&&$1$&&&1&&\cr
\hline
$1$ &&&1&&& &&&& &&\cr
\hline
$2$&&&&&& &&&$\bar 1$&& &$\bar 1$\cr
\hline
$3$&1&&&& &&&&& &&\cr
\hline
\end{tabular}
$$}

\bigskip\noindent{\small
D.L. was partly supported by l'IHES and MPIMiS.}

\end{document}